\documentclass[preprint,12pt]{elsarticle}

\usepackage{amssymb}
\usepackage{amsthm}
\usepackage{amsmath}
\usepackage{xcolor}
\newtheorem{theorem}{Theorem}[section]
\newtheorem{proposition}{Proposition}[section]
\newtheorem{example}{Example}[section]%
\newtheorem{remark}{Remark}[section]%
\newtheorem{lemma}{Lemma}[section]%
\begin{document}
	
	\begin{frontmatter}
		
		
		
		\title{On constacyclic codes over a class of non-chain rings}
		
		
		\author{Nikita Jain\fnref{label1}}
		\ead{nikitajain.phd19appsc@pec.edu.in}
		
		\author{Sucheta Dutt\corref{cor1}\fnref{label2}}
		\cortext[cor1]{Corresponding author}
		\ead{sucheta@pec.edu.in}
		
		\author{Ranjeet Sehmi\fnref{label3}}
		\ead{rsehmi@pec.edu.in}
		
		\affiliation[label1]{organization={Punjab Engineering College (Deemed to be University), Chandigarh},
			addressline={Sector 12, Chandigarh}, 
			city={Chandigarh},
			postcode={160012}, 
			state={Punjab},
			country={India}}
		\affiliation[label2]{organization={Punjab Engineering College (Deemed to be University), Chandigarh},
			addressline={Sector 12, Chandigarh}, 
			city={Chandigarh},
			postcode={160012}, 
			state={Punjab},
			country={India}}
		\affiliation[label3]{organization={Punjab Engineering College (Deemed to be University), Chandigarh},
			addressline={Sector 12, Chandigarh}, 
			city={Chandigarh},
			postcode={160012}, 
			state={Punjab},
			country={India}}
		
		\begin{abstract}
		In this paper, a unique form of generators of a constacyclic code of arbitrary length over a non-chain ring of the type $\mathtt{R_{_{\theta}}}=Z_{4}+\nu Z_{4}, \nu^{2}=\theta \in Z_{4}+\nu Z_{4}$ has been obtained. Further, rank and cardinality of a constacyclic code of arbitrary length over a non-chain ring of the type $\mathtt{R_{_{\theta}}}$ have been obtained by determining a minimal spanning set of the code.  Also, necessary and sufficient conditions for a constacyclic code of arbitrary length over a non-chain ring of the type $\mathtt{R_{_{\theta}}}$ to be reversible have been determined. Examples have also been presented in support of our results.
			
		\end{abstract}

		
	\end{frontmatter}
	
	
	\section{Introduction}
	\label{}
Constacyclic codes were first introduced by Berlekamp \cite{Berlekamp} over fields. Constacyclic codes of various lengths over different rings have been studied by many researchers with diverse approaches. For reference, see (\cite{CCC ref 6,CCC ref2,CCC ref1,CCC ref 5,CCC ref 4,CCC ref 2,CCC ref 3,lemma 3ps hamming distance,CCC ref 10,CCC ref 8,CCC ref 9,CCC ref 7}). The generators of $(1+u)$- constacyclic codes of odd length over $F_{2}+uF_{2}$, $u^{2}=0$ have been obtained by Qian et al. \cite{CC1}. The generators of $(1+u)$- constacyclic codes of arbitrary length over the same ring $F_2 +uF_2, u^{2}=0$ have been obtained by Abualrub et al. \cite{CC2}. The structure of $(1+v)$- constacyclic codes of odd length over $F_2 +uF_2 + vF_2 +uvF_2$, where $ u^{2}= v^{2}=0$ have been determined by Karadeniz et al. \cite{CC3}.
The generators of $(1+2\nu)$- constacyclic codes of odd length over $Z_{4}+\nu Z_{4}$ for $\nu^{2}=0$ have been obtained by M. Ashraf and G. Mohammad \cite{chapter 5(1)}. 
Further, the generators of $(1+\nu)$- constacyclic codes of arbitrary length over $Z_{4}+\nu Z_{4}$ for $\nu^{2}=0$ have been determined by Haeifeng Yu et al. \cite{chapter 5(2)}. 
The generators of $(2+\nu)$- constacyclic codes of odd length over $Z_{4}+\nu Z_{4}$ for $\nu^{2}=1$ have been determined by $\ddot{O}$zen et al. \cite{cyclic and some consta cyclic codes by ozen}. 
The generators of
$(1+2\nu),(3+2\nu)$ and $(2+3\nu)$- constacyclic codes of odd length over $Z_{4}+\nu Z_{4}$ for $\nu^{2}=1$ have been obtained in \cite{on constacyclic,on some constacyclic}. Further, the generators of $(1+2\nu)$ and $(3+2\nu)$- constacyclic codes of odd length over $Z_{4}+\nu Z_{4},\nu^{2} =\nu$ have been determined by H. Dinh et al. \cite{Dinh CC}.

Most of the work done on constacyclic codes over a non-chain ring of the type $\mathtt{R_{_{\theta}}}$ deals with constacyclic codes of odd length. Moreover, all available studies as listed above deal with a $\lambda$- constacyclic code for a particular value of $\lambda.$ With an aim to address these reasearch gaps, in this paper, a unique form of generators of a constacyclic code of arbitrary length along with its rank and cardinality over a non-chain ring of the type $\mathtt{R_{_{\theta}}}$ has been obtained. Also, necessary and sufficient conditions for a constacyclic code of arbitrary length over a non-chain ring of the type $\mathtt{R_{_{\theta}}}$ to be reversible have been determined.	

	\section{Preliminaries}
	\label{}
	Let $\mathtt{R}$ be a finite commutative ring with unity. For a unit $\lambda$ in $\mathtt{R}$ and an element $\mathtt{s}=(\mathtt{s_{0}},\mathtt{s_{1}},\cdots,\mathtt{s_{n-1}})$ in $\mathtt{C}$, the map $\tau_{\lambda}(\mathtt{s_{0}},\mathtt{s_{1}},\cdots,\mathtt{s_{n-1}})=(\lambda \mathtt{s_{n-1}},\mathtt{s_{0}},\mathtt{s_{1}},\\\cdots,\mathtt{s_{n-2}})$ is called the $\lambda$-constacyclic shift of $\mathtt{s}.$ If $\mathtt{C}$ is closed under the $\lambda$-constacyclic shift $\tau_{\lambda}$, i.e, $\tau_{\lambda} (\mathtt{C}) = \mathtt{C}$, then the code $\mathtt{C}$ is called a $\lambda$-constacyclic code of length $n$ over $\mathtt{R}$. For $\lambda=1,$ $\mathtt{C}$ is said to be a cyclic code and for $\lambda=-1,$ $\mathtt{C}$ is said to be a negacyclic code of length $n$ over $\mathtt{R}.$ It is well known that the ideals of $\mathtt{R}[\mathrm{z}]/\langle \mathrm{z}^{n}-\lambda\rangle $ and $\lambda$-constacyclic codes of length $n$ over $\mathtt{R}$ are in one to one correspondence. A finite commutative ring $\mathtt{R}$ is called a chain ring if all its ideals form a chain under the inclusion relation. Otherwise, $\mathtt{R}$ is called a non-chain ring. A set $\mathtt{S}$ of elements of a constacyclic code $\mathtt{C}$ over a finite commutative ring $\mathtt{R}$ is called a spanning set of $\mathtt{C}$ if each element of $\mathtt{C}$ can be written as a linear combination of elements of $\mathtt{S}$ with coefficients in $\mathtt{R}.$ A spanning set $\mathtt{S}$ of a constacyclic code $\mathtt{C}$ over a finite commutative ring $\mathtt{R}$ is called a minimal spanning set of $\mathtt{C}$ if no proper subset of $\mathtt{S}$ spans $\mathtt{C}.$ It can be easily proved that the number of elements in any two minimal spanning sets of $\mathtt{C}$ remains the same. The rank of a constacyclic code $\mathtt{C}$ over a finite commutative ring $\mathtt{R}$ is defined as the number of elements in a minimal spanning set of $\mathtt{C}.$ The rank of $\mathtt{C}$ is denoted by $rank(\mathtt{C}).$ The cardinality of $\mathtt{C},$ denoted by $\vert \mathtt{C} \vert,$ is defined as the total number of elements of $\mathtt{C}.$  \\
	
	The rings $Z_{4}+\nu Z_{4}$, $\nu^{2} \in Z_{4}+\nu Z_{4}$ have been classified into chain rings and non-chain rings by Adel Alahmadi et al.\cite{classification base paper}. They have proved that $Z_{4}+\nu Z_{4}$ is a chain ring for $\nu^{2} \in \{2,3,1+\nu,1+2\nu,1+3\nu,2+2\nu,3+\nu,3+3\nu\}$ and is a non-chain ring for $\nu^{2} \in \{0,1,\nu,2\nu,3\nu,2+\nu,2+3\nu,3+2\nu\}.$
	It is easy to see that among the eight non-chain rings $Z_{4}+\nu Z_{4},\nu^{2} \in \{0,1,\nu,2\nu,3\nu,2+\nu,2+3\nu,3+2\nu\}$  some exhibit isomorphism with each other. The rings $Z_{4}+\nu Z_{4}$, $\nu^{2}=0$ and $Z_{4}+\nu Z_{4}$, $\nu^{2}=3+2\nu$ are found to be isomorphic. Additionally, the ring $Z_{4}+\nu Z_{4}$, $\nu^{2}=\nu$ is isomorphic to the rings $Z_{4}+\nu Z_{4}$, where $\nu^{2} \in \{3\nu,2+\nu,2+3\nu\}.$ Similarly, the ring $Z_{4}+\nu Z_{4}$ with $\nu^{2}=1$ is isomorphic to the ring $Z_{4}+\nu Z_{4}$ with $\nu^{2}=2\nu.$
	
	Thus to prove any results for codes over all the non-chain rings of the type $Z_{4}+\nu Z_{4},$ where $\nu^{2} \in \{0,1,\nu,2\nu,3\nu,2+\nu,2+3\nu,3+2\nu\}$ it is enough to prove the results for codes over the rings $Z_{4}+\nu Z_{4},$ where $\nu^{2} \in \{0,1,\nu\}.$\\
	
	Let $\mathtt{R}_{_{\theta}}=Z_{4}+\nu Z_{4},$ where $\nu^{2}=\theta \in \{0,1,\nu\}.$ 
	Define\\
	\[ k_{_{\theta}}= \begin{cases}
		\nu & ;~\theta \in \{0,\nu\}\\
		1+\nu & ;~\theta =1\\
	\end{cases} \]
	
	\begin{lemma}{\label{iso homo lemma}}\cite{self paper}
		The map $\phi_{_{\theta}}: \mathtt{R}_{_{\theta}} \longrightarrow Z_{4}$ defined as $\phi_{_{\theta}}(x) =x~ (mod~ k_{_{\theta}})$ is a ring homomorphism for $\theta \in \{0,1,\nu\}.$
	\end{lemma}
	Let $u_{_{\theta}}$ be a unit element in the ring $\mathtt{R_{_{\theta}}}; \theta \in \{0,1,\nu\}.$ 
	Let $\mathtt{C_{_{\theta}}}$ be a $u_{_{\theta}}$- constacyclic code of length $n$ over $\mathtt{R_{_{\theta}}}, \theta \in \{0,1,\nu\}.$
	The map $\phi_{_{\theta}}$ can be naturally extended to a map $\phi_{_{\theta}}$: $\mathtt{C}_{_{\theta}} \rightarrow Z_{4}[z]/\langle z^{n}-u_{_{\theta}}(mod ~k_{_{\theta}}) \rangle$ as $\phi_{_{\theta}}\big(\mathtt{s}_{_{0}}+\mathtt{s}_{_{1}}z+\cdots+\mathtt{s}_{_{n-1}}z^{n-1}\big) = \phi_{_{\theta}}(\mathtt{s}_{_{0}})+\phi_{_{\theta}}(\mathtt{s}_{_{1}})z+\cdots+\phi_{_{\theta}}(\mathtt{s}_{_{n-1}})z^{n-1}.$ Then, kernel of $\phi_{_{\theta}}$ is defined as
	$ker_{_{\theta}}=\bigl\{a(z) \in \mathtt{C}_{_{\theta}}$ such that $\phi_{_{\theta}}(a(z))=0\bigr\}.$ Further, Torsion of $\mathtt{C_{_{\theta}}}$ is defined as Tor$(\mathtt{C}_{_{\theta}})= \Bigl\{a(z)\in \dfrac{Z_{4}[z]}{\langle z^{n}-u_{_{\theta}}(mod ~k_{_{\theta}})\rangle}: k_{_{\theta}} a(z)\in \mathtt{C}_{_{\theta}}\Bigr\}$ and Residue of $\mathtt{C}_{_{\theta}}$ is defined as
	Res$(\mathtt{C}_{_{\theta}})= \Bigl\{a(z)\in \dfrac{Z_{4}[z]}{\langle z^{n}-u_{_{\theta}}(mod ~k_{_{\theta}})\rangle}:a(z)+k_{_{\theta}} b(z) \in \mathtt{C}_{_{\theta}}$ for some $b(z) \in \dfrac{Z_{4}[z]}{\langle z^{n}-u_{_{\theta}}(mod ~k_{_{\theta}})\rangle}\Bigr\}.$ Clearly, Tor$(\mathtt{C}_{_{\theta}})$ and Res$(\mathtt{C}_{_{\theta}})$ are ideals of the ring $\dfrac{Z_{4}[z]}{\langle z^{n}-u_{_{\theta}}(mod ~k_{_{\theta}})\rangle}.$ Define the following sets
	
	$\mathtt{C}_{_{\theta_{1}}} = $ Res $($Res $(\mathtt{C}_{_{\theta}})) = \mathtt{C}_{_{\theta}}$ mod $(2,k_{_{\theta}}),$
	
	$\mathtt{C}_{_{\theta_{2}}} = $ Tor $($Res $(\mathtt{C}_{_{\theta}})) = \{a(z)\in Z_{2}[z]:2a(z) \in \mathtt{C}_{_{\theta}}$ mod $k_{_{\theta}}\},$
	
	$\mathtt{C}_{_{\theta_{3}}} = $ Res $($Tor $(\mathtt{C}_{_{\theta}}))= \{a(z)\in Z_{2}[z]:k_{_{\theta}}a(z) \in \mathtt{C}_{_{\theta}}$ mod $2k_{_{\theta}}\},$
	
	$\mathtt{C}_{_{\theta_{4}}} = $ Tor $($Tor $(\mathtt{C}_{_{\theta}}))= \{a(z)\in Z_{2}[z]:2k_{_{\theta}}a(z) \in \mathtt{C}_{_{\theta}}\}.$
	
	It can be easily seen that the sets $\mathtt{C}_{_{\theta_{1}}},\mathtt{C}_{_{\theta_{2}}},\mathtt{C}_{_{\theta_{3}}},\mathtt{C}_{_{\theta_{4}}}$ defined above are ideals  of the ring $\dfrac{Z_{2}[z]}{\langle z^{n}-1 \rangle}$ and are therefore principally generated.
	
	\section{A unique form of generators of a constacyclic code over $\mathtt{R_{_{\theta}}}$}\label{sec2}
	In this section, a unique form of generators of a constacyclic code of arbitrary length over $\mathtt{R_{_{\theta}}}, \theta \in \{0,1,\nu\}$ has been obtained. It can be easily observed that if $\mathtt{C_{_{u_{_{\theta}}}}}= \mathtt{C_{_{1,u_{_{\theta}}}}}+k{_{_{\theta}}}\mathtt{C_{_{2,u_{_{\theta}}}}}$ is a $u_{_{\theta}}$- constacyclic code of arbitrary length $n$ over the ring $\mathtt{R_{_{\theta}}}$ for $\theta \in \{0,1,\nu\},$ then $\mathtt{C_{_{1,u_{_{\theta}}}}}$ is either a cyclic code or a negacyclic code of length $n$ over $Z_{4}.$ Depending on whether $\mathtt{C_{_{1,u_{_{\theta}}}}}$ is a cyclic code or a negacyclic code over $Z_{4}$, we categorize all the unit elements $u_{_{\theta}}$ of $\mathtt{R_{_{\theta}}}$ into two groups $u_{_{\theta}}=\alpha_{_{\theta}}$ and $u_{_{\theta}}=\beta_{_{\theta}}$ as mentioned below in Table 2.1.\\
	
	\begin{center}
		\begin{tabular}{|c|c|c|}
			\hline
			$\theta$ & $u_{_{\theta}}=\alpha_{_{\theta}}$ & $u_{_{\theta}}=\beta_{_{\theta}}$ \\
			\hline
			$0$ & $1+\nu,1+2\nu,1+3\nu$ & $3,3+2\nu,3+3\nu,3+\nu$\\
			\hline
			$1$ & $2+\nu,3+2\nu,3\nu$ & $3,\nu,1+2\nu,2+3\nu$\\
			\hline
			$\nu$ & $1+2\nu$ & $3+2\nu$\\
			\hline
		\end{tabular}	
	\end{center}
	\begin{center}
		Table 2.1
	\end{center}
	It is noted that for all values of $u_{_{\theta}}=\alpha_{_{\theta}},$ $\mathtt{C_{_{1,u_{_{\theta}}}}}$ is a cyclic code and for all values of $u_{_{\theta}}=\beta_{_{\theta}},$ $\mathtt{C_{_{1,u_{_{\theta}}}}}$ is a negacyclic code over $Z_{4}.$
	
	We shall require the structure of a cyclic code of arbitrary length over $Z_{4}$ given by Abualrub and Siap \cite{gen over Z4 abul} to prove our results. For the sake of clarity and completeness, we recall the relevant result from \cite{gen over Z4 abul} as Lemma \ref{lemma gen over $Z_{4}$} below. 
	
	\begin{lemma} \cite{gen over Z4 abul} \label{lemma gen over $Z_{4}$} Let $\mathtt{C}$ be a cyclic code of arbitrary length $n$ over $Z_{4}$. Then $\mathtt{C}=\langle g(z)+2p(z),2a(z)\rangle,$ where $ g(z),a(z)$ and $p(z)$ are binary polynomials such that $a(z)\vert g(z)\vert z^{n}-1$ and either $p(z)=0$ or $a(z)\vert p(z)\frac{z^{n}-1}{g(z)}$ with deg $a(z)>$ deg $p(z)$.
	\end{lemma}
	
	A unique form of generators of a cyclic code of arbitrary length over a non-chain ring of the type $Z_{4}+\nu Z_{4},$ where $\nu^{2} \in Z_{4}+\nu Z_{4}$ has been obtained by Jain et al. \cite{self paper}. It is observed that a similar approach can be used to obtain a unique form of generators of an $\alpha_{_{\theta}}$- constacyclic code of arbitrary length over a non-chain ring of the type $Z_{4}+\nu Z_{4},$ where $\nu^{2} \in Z_{4}+\nu Z_{4}.$
	
	In the following theorems, a unique form of generators of an $\alpha_{_{\theta}}$- constacyclic code $\mathtt{C}_{_{\theta}}$ of arbitrary length $n$ over the ring $\mathtt{R}_{_{\theta}}, \theta \in \{0,1,\nu\}$ has been established.
	\begin{theorem}\label{consta from cyclic}
		Let $\mathtt{C_{_{\theta}}}$ be an $\alpha_{_{\theta}}$- constacyclic code of arbitrary length $n$ over $\mathtt{R_{_{\theta}}}, \theta \in \{0,1,\nu\}.$ Then $\mathtt{C_{_{\theta}}} =\langle h_{_{\theta_{1}}}(z),h_{_{\theta_{2}}}(z),h_{_{\theta_{3}}}(z),h_{_{\theta_{4}}}(z)\rangle,$ where $h_{_{\theta_{i}}}(z) \in \dfrac{\mathtt{R_{_{\theta}}}[z]}{\langle z^{n}-\alpha_{_{\theta}} \rangle}$ for $1 \leq i \leq 4$. Further, $h_{_{\theta_{1}}}(z)= h_{_{11}}(z)+2h_{_{12}}(z)+k_{_{\theta}}h_{_{13}}(z)+2k_{_{\theta}}h_{_{14}}(z)$, $h_{_{\theta_{2}}}(z)=2h_{_{22}}(z)+k_{_{\theta}}h_{_{23}}(z)+2k_{_{\theta}}h_{_{24}}(z)$, $h_{_{\theta_{3}}}(z)=k_{_{\theta}}h_{_{33}}(z)+2k_{_{\theta}}h_{_{34}}(z)$ and $h_{_{\theta_{4}}}(z)=2k_{_{\theta}}h_{_{44}}(z)$ such that the polynomials $h_{_{ij}}(z)$ are in $Z_{2}[z]/{\langle z^{n}-1\rangle}$ for $1 \leq i \leq j \leq 4.$ Further,
		\begin{equation}\label{eqn1cc}
			h_{_{22}}(z)\big\vert h_{_{11}}(z)\big\vert z^{n}-1,   
		\end{equation}
		\begin{equation}
			either~h_{_{12}}(z)=0 ~or~ h_{_{22}}(z)\big\vert h_{_{12}}(z)\frac{z^{n}-1}{h_{_{11}}(z)}~~ with~deg ~h_{_{22}}(z)>  deg ~h_{_{12}}(z), 
		\end{equation}
		\begin{equation}\label{equation 2.3cc}
			h_{_{44}}(z) \big\vert h_{_{33}}(z)\big\vert z^{n}-1,
		\end{equation}
		\begin{equation}\label{eqn4cc}
			either~h_{_{34}}(z)=0 ~or~h_{_{44}}(z)\big\vert h_{_{34}}(z)\frac{z^{n}-1}{h_{_{33}}(z)}  ~~with~deg ~h_{_{44}}(z)>  deg ~h_{_{34}}(z). 
		\end{equation}
	\end{theorem}
	\begin{proof} 	Let $\mathtt{C}_{_{\theta}}$ be an $\alpha_{_{\theta}}$- constacyclic code of length $n$ over $\mathtt{R}_{_{\theta}}$, $\theta \in \{0,1,\nu\}.$ Let $\phi_{_{\theta}}$ be the ring homomorphism as defined in Lemma \ref{iso homo lemma} for $\theta \in \{0,1,\nu\}.$ 
		Restrict $\phi_{_{\theta}}$ to $\mathtt{C}_{_{\theta}}.$ It is easy to see that, $\phi_{_{\theta}}(\mathtt{C}_{_{\theta}})$ is a cyclic code of length $n$ over $Z_{4}.$  Using Lemma \ref{lemma gen over $Z_{4}$}, we get $\phi_{_{\theta}}(\mathtt{C}_{_{\theta}})$ = $\langle h_{_{11}}(z)+2h_{_{12}}(z),2h_{_{22}}(z)\rangle$, where $h_{_{22}}(z)\big\vert h_{_{11}}(z)\big\vert z^{n}-1$ and either $h_{_{12}}(z)=0$ or $h_{_{22}}(z)\big\vert h_{_{12}}(z)\frac{z^{n}-1}{h_{_{11}}(z)}$ with deg $h_{_{22}}(z)>$ deg $h_{_{12}}(z).$  
		Also $ker_{_{\theta}}=\langle k_{_{\theta}} \rangle.$ Consider $I=\biggl\{a(z)\in \dfrac{Z_{4}[z]}{\langle z^{n}-1\rangle}: k_{_{\theta}}a(z)\in ker_{_{\theta}}\biggr\}.$ Clearly, $I$ is an ideal of $\dfrac{Z_{4}[z]}{\langle z^{n}-1\rangle}.$ Thus, $ker_{_{\theta}}$ is an ideal of $k_{_{\theta}}Z_{4}[z]/\langle z^{
			n}-1\rangle.$ Therefore, $ker_{_{\theta}}$ is of the form $k_{_{\theta}}J,$ where $J$ is a cyclic code of length $n$ over $Z_{4}.$ 
		Therefore, using Lemma \ref{lemma gen over $Z_{4}$}, we obtain $ker_{_{\theta}}$ = $k_{_{\theta}}\langle h_{_{33}}(z)+2h_{_{34}}(z),2h_{_{44}}(z)\rangle$, where  $h_{_{44}}(z)\big\vert h_{_{33}}(z)\big\vert z^{n}-1$ and either $h_{_{34}}(z)=0$ or $h_{_{44}}(z)\big\vert h_{_{34}}(z)\frac{z^{n}-1}{h_{_{33}}(z)}$ with deg $h_{_{44}}(z)>$ deg $h_{_{34}}(z).$
		
		It follows that
		$\mathtt{C}_{_{\theta}}$ is generated by $h_{_{\theta_{1}}}(z)= h_{_{11}}(z)+2h_{_{12}}(z)+k_{_{\theta}}h_{_{13}}(z)+2k_{_{\theta}}h_{_{14}}(z)$, $h_{_{\theta_{2}}}(z)=2h_{_{22}}(z)+k_{_{\theta}}h_{_{23}}(z)+2k_{_{\theta}}h_{_{24}}(z)$, $h_{_{\theta_{3}}}(z)=k_{_{\theta}}h_{_{33}}(z)+2k_{_{\theta}}h_{_{34}}(z)$ and $h_{_{\theta_{4}}}(z)=2k_{_{\theta}}h_{_{44}}(z)$ such that $h_{_{\theta_{i}}}(z) \in \dfrac{\mathtt{R_{_{\theta}}}[z]}{\langle z^{n}-\alpha_{_{\theta}} \rangle}$ and the polynomials $h_{_{ij}}(z)$ are in $Z_{2}[z]/{\langle z^{n}-1\rangle}$ for $1 \leq i \leq j \leq 4$ and satisfy the conditions (\ref{eqn1cc})-(\ref{eqn4cc}).
	\end{proof}
	\begin{remark}
		It can be easily seen that for an $\alpha_{_{\theta}}$- constacyclic code $\mathtt{C_{_{\theta}}}$ of arbitrary length $n$ over the ring $\mathtt{R_{_{\theta}}}, \theta \in \{0,1,\nu\}$; the cyclic codes $\mathtt{C}_{_{\theta_{1}}},\mathtt{C}_{_{\theta_{2}}},\mathtt{C}_{_{\theta_{3}}}$ and $\mathtt{C}_{_{\theta_{4}}}$ are generated respectively by the polynomials $h_{_{11}}(z),h_{_{22}}(z),h_{_{33}}(z)$ and $h_{_{44}}(z).$ Thus $h_{_{11}}(z),h_{_{22}}(z),h_{_{33}}(z)$ and $h_{_{44}}(z)$ can be considered to be the unique minimal degree polynomial generators of the principally generated codes $\mathtt{C}_{_{\theta_{1}}},\mathtt{C}_{_{\theta_{2}}},\mathtt{C}_{_{\theta_{3}}}$ and $\mathtt{C}_{_{\theta_{4}}}$ respectively.   
	\end{remark}
		\begin{theorem} \label{theorem str 2 cc} Let $\mathtt{C}_{_{\theta}} = \langle h_{_{\theta_{1}}}(z),h_{_{\theta_{2}}}(z),h_{_{\theta_{3}}}(z),h_{_{\theta_{4}}}(z)\rangle$ be an $\alpha_{_{\theta}}$- constacyclic code of arbitrary length $n$ over the ring $\mathtt{R}_{_{\theta}}, \theta \in \{0,1,\nu\},$ where  $h_{_{\theta_{i}}}(z) \in \dfrac{\mathtt{R_{_{\theta}}}[z]}{\langle z^{n}-\alpha_{_{\theta}} \rangle}$ for $1 \leq i \leq 4$ are polynomials as defined in Theorem \ref{consta from cyclic}. Then there exists a set of generators  $\{t_{_{\theta_{1}}}(z),t_{_{\theta_{2}}}(z),t_{_{\theta_{3}}}(z),t_{_{\theta_{4}}}(z)\}$ of $\mathtt{C}_{_{_\theta}},$ where $t_{_{\theta_{1}}}(z)\\= t_{_{11}}(z)+2t_{_{12}}(z)+k_{_{\theta}}t_{_{13}}(z)+2k_{_{\theta}}t_{_{14}}(z)$, $t_{_{\theta_{2}}}(z)=2t_{_{22}}(z)+k_{_{\theta}}t_{_{23}}(z)+2k_{_{\theta}}t_{_{24}}(z)$, $t_{_{\theta_{3}}}(z)=k_{_{\theta}}t_{_{33}}(z)+2k_{_{\theta}}t_{_{34}}(z)$ and $t_{_{\theta_{4}}}(z)=2k_{_{\theta}}t_{_{44}}(z)$ such that the polynomials $t_{_{ij}}(z)$ for $1 \leq i \leq j \leq 4$ are in $Z_{2}[z]/{\langle z^{n}-1 \rangle}$ and satisfy the conditions (\ref{eqn1cc})-(\ref{eqn4cc}) as defined in Theorem \ref{consta from cyclic}. Further $t_{_{ii}}(z)$ are unique minimal degree polynomial generators of $\mathtt{C}_{_{\theta_{i}}}, 1 \leq i \leq 4.$ Also, either $t_{_{ij}}(z)=0$ or deg $t_{_{ij}}(z)<$ deg $t_{_{jj}}(z)$ for $ 1\leq i \leq 3, i<j\leq 4.$ 
	\end{theorem}
	\begin{proof}
		The proof follows a similar approach to the proof of Theorem 4 \cite{self paper}.
	\end{proof}
	\begin{theorem} \label{theorem str 3cc}	Let $\mathtt{C}_{_{\theta}}=\langle t_{_{\theta_{1}}}(z),t_{_{\theta_{2}}}(z),t_{_{\theta_{3}}}(z),t_{_{\theta_{4}}}(z)\rangle$ be an $\alpha_{_{\theta}}$- constacyclic code of arbitrary length $n$ over the ring $\mathtt{R}_{_{\theta}}, \theta \in \{0,1,\nu\}$, where $t_{_{\theta_{i}}}(z) \in \dfrac{\mathtt{R_{_{\theta}}}[z]}{\langle z^{n}-\alpha_{_{\theta}} \rangle}$ for $1 \leq i \leq 4$ are as obtained in Theorem \ref{theorem str 2 cc} above. 
		Then the polynomials $t_{_{\theta_{i}}}(z)$ for $1 \leq i \leq 4$ are uniquely determined.
	\end{theorem}
	\begin{proof}
		The proof follows the same lines as the proof of Theorem 5 \cite{self paper}.
	\end{proof}
	
	The following theorem which gives some divisibility properties of polynomials $t_{_{ij}}(z), 1 \leq i \leq j \leq 4$ in $Z_{2}[z]/\langle z^{n}-1\rangle$ can be proved easily through simple calculations. 
	\begin{theorem} \label{theorem str div cond cc} Let $\mathtt{C}_{_{\theta}}=\langle t_{_{\theta_{1}}}(z),t_{_{\theta_{2}}}(z),t_{_{\theta_{3}}}(z),t_{_{\theta_{4}}}(z)\rangle$ be an $\alpha_{_{\theta}}$- constacyclic code of arbitrary length $n$ over the ring $\mathtt{R}_{_{\theta}}, \theta \in \{0,1,\nu\},$ where the generators $t_{_{\theta_{1}}}(z)= t_{_{11}}(z)+2t_{_{12}}(z)+k_{_{\theta}}t_{_{13}}(z)+2k_{_{\theta}}t_{_{14}}(z)$, $t_{_{\theta_{2}}}(z)=2t_{_{22}}(z)+k_{_{\theta}}t_{_{23}}(z)+2k_{_{\theta}}t_{_{24}}(z)$, $t_{_{\theta_{3}}}(z)=k_{_{\theta}}t_{_{33}}(z)+2k_{_{\theta}}t_{_{34}}(z)$ and $t_{_{\theta_{4}}}(z)=2k_{_{\theta}}t_{_{44}}(z)$ are in the unique form as in Theorem \ref{theorem str 3cc}. Then the following divisibility relations hold over the ring $Z_{2}.$
		\begin{itemize}
			\item [(i)] $t_{_{33}}(z)\Big\vert\frac{z^{n}-1}{t_{_{11}}(z)}\Big(t_{_{13}}(z)-\frac{t_{_{12}}(z)}{t_{_{22}}(z)}t_{_{23}}(z)\Big)$
			\item [(ii)] $t_{_{44}}(z)\big\vert t_{_{23}}(z)$
			\item [(iii)] $t_{_{33}}(z)\Big \vert\frac{t_{_{11}}(z)}{t_{_{22}}(z)}t_{_{23}}(z)$
			\item [(iv)] $t_{_{44}}(z)\Big \vert \frac{z^{n}-1}{t_{_{22}}(z)}\Big(t_{_{24}}(z)-\frac{t_{_{23}}(z)}{t_{_{33}}(z)}t_{_{34}}(z)\Big)$
			\item [(v)] 	$t_{_{44}}(z)\Big \vert t_{_{13}}(z)-\frac{t_{_{11}}(z)}{t_{_{22}}(z)}t_{_{24}}(z)+\frac{t_{_{11}}(z)}{t_{_{22}}(z)t_{_{33}}(z)}t_{_{23}}(z)t_{_{34}}(z)$
			\item [(vi)] $t_{_{44}}(z)\Big \vert\frac{z^{n}-1}{t_{_{11}}(z)}\Big(t_{_{14}}(z)-\frac{t_{_{12}}(z)}{t_{_{22}}(z)}t_{_{24}}(z)+\frac{-t_{_{13}}(z)+\frac{t_{_{12}}(z)t_{_{23}}(z)}{t_{_{22}}(z)}}{t_{_{33}}(z)}t_{_{34}}(z)\Big)$
			\item[(vii)] $t_{_{44}}(z)\big \vert t_{_{11}}(z),$ 	$t_{_{44}}(z)\big \vert t_{_{22}}(z),$ 	$t_{_{44}}(z)\big \vert t_{_{33}}(z),$ 	$t_{_{33}}(z)\big \vert t_{_{11}}(z)$ for $\theta=0$
			\item [(viii)] 	
			$t_{_{33}}(z)\big \vert t_{_{11}}(z),$
			$t_{_{44}}(z)\big \vert t_{_{11}}(z),$
			$t_{_{44}}(z)\big \vert t_{_{22}}(z)+t_{_{23}}(z)$ for $\theta =1$
			\item [(ix)] $g_{_{44}}(z)\big \vert t_{_{12}}(z)+t_{_{13}}(z)-\frac{t_{_{11}}(z)}{t_{_{33}}(z)}t_{_{34}}(z)$ for $\theta =1$
			\item [(x)]	$t_{_{44}}(z)\big \vert t_{_{13}}(z), 	t_{_{44}}(z)\big \vert t_{_{11}}(z)$  for $\theta =\nu.$
		\end{itemize}
	\end{theorem}
	The algebraic structure of all $\lambda$- constacyclic codes of length $np^{s}$ over a finite commutative chain ring $R_{e}$ with nilpotency index $e\geq 2$ has been given by Sharma and Sidana \cite{anuradha}. We recall below the relevant results from \cite{anuradha} which we shall use to establish the algebraic structure of all $\beta_{_{\theta}}$- constacyclic codes over $\mathtt{R_{_{\theta}}}, \theta \in \{0,1,\nu\}.$ 
	
	\begin{proposition}\label{prop techmuller set}(Proposition 2.2,\cite{anuradha}) Let $R_{e}$ be a finite commutative chain ring with unity 1 and with the maximal ideal as $M= 
		\langle \gamma \rangle,$ where $e \geq 2$ is the nilpotency index of the generator $\gamma$ of $M$. Next, let $\overline{R_{e}}=R_{e}/ \langle \gamma \rangle$ be the residue field of $R_{e}.$ As $\overline{R_{e}}$ is a finite field, $char \overline{R_{e}}$ is a prime number, say $p.$ Let us suppose that $\vert \overline{R_{e}}\vert=p^{m}$ for some positive integer $m.$ Then we have the following:
		\begin{itemize}
			\item [(i)] The characteristic of $R_{e}$ is $p^{a},$ where $ 1 \leq a \leq e.$ Moreover, $\vert R_{e}\vert =\vert \overline{R_{e}}\vert ^{e}=p^{me}.$
			\item [(ii)] There exists an element $\zeta \in R_{e}$ having the multiplicative order $p^{m}-1.$ Moreover, each element $r \in R_{e}$ can be uniquely expressed as $r=r_{0}+r_{1}\gamma+\cdots+r_{e-1}\gamma^{e-1},$ where $r_{i} \in T_{e}=\{0,1,\zeta,\cdots,\zeta^{p^{m}-2}\}$ for $0 \leq i \leq e-1.$ (The set $T_{e}=\{0,1,\zeta,\cdots,\zeta^{p^{m}-2}\}$ is called the Teichm$\ddot{u}$ller set of $R_{e}).$
			\item [(iii)] Let $r=r_{0}+r_{1}\gamma+\cdots+r_{e-1}\gamma^{e-1},$ where $r_{i} \in T_{e}=\{0,1,\zeta,\cdots,\zeta^{p^{m}-2}\}$ for $0 \leq i \leq e-1.$ Then $r$ is a unit in $R_{e}$ if and only if $r_{0}\neq 0$. Moreover, there exists $\alpha_{0} \in T_{e}$ satisfying $\alpha_{0}^{p^{s}}=r_{0}.$
		\end{itemize}
	\end{proposition}
	
	\begin{lemma}(Theorem 4.1(c),\cite{anuradha})\label{anuradha nega cyclic result}
		Let $R_{e}$ be a finite commutative chain ring with unity $1$ and with the maximal ideal as $M= 
		\langle \gamma \rangle,$ where $e \geq 2$ is the nilpotency index of the generator $\gamma$ of $M$. Next, let $\overline{R_{e}}=R_{e}/ \langle \gamma \rangle$ be the residue field of $R_{e}.$ As $\overline{R_{e}}$ is a finite field, $char \overline{R_{e}}$ is a prime number, say $p.$ Let $\lambda= \theta+\gamma \omega,$ where $\theta \in T_{e}$  and $\omega $ is unit in $R_{e}.$ As $\theta \in T_{e},$ by Proposition \ref{prop techmuller set}(iii), there exists $\lambda_{_{0}} \in T_{e} $ satisfying $\theta =\lambda_{_{0}}^{p^{s}}.$ Then all the $\lambda$- constacyclic codes over $R_{e}$ of length $np^{s},$ where $s$ and $n$ are positive integers such that $gcd(n,p)=1$ are given by $\langle (x^{n}-\lambda_{_{0}})^\nu\rangle,$ where $ 0 \leq \nu \leq ep^{s}.$
	\end{lemma}
	Taking $R_{e}$ as $Z_{4}$ and $\lambda=-1$ in Lemma \ref{anuradha nega cyclic result} above, we obtain the generators of all negacyclic codes of length $n.2^{s}$ over $Z_{4}.$ We state this special case of Lemma \ref{anuradha nega cyclic result} as Lemma \ref{nega cyclic gen over $Z_{4}$} below.
	
	\begin{lemma}\label{nega cyclic gen over $Z_{4}$}
		All the negacyclic codes of length $n.2^{s}$ over $Z_{4}$ are given by $\langle (z^{n}-1)^{t} \rangle, 0 \leq t \leq 2^{s+1},$ where $s$ and $n$ are positive integers and $n$ is odd. 
	\end{lemma}
	In the following theorem, the generators of a $\beta_{_{\theta}}$- constacyclic code of arbitrary length $n$ over $\mathtt{R_{_{\theta}}}, \theta \in \{0,1,\nu\}$ have been obtained. 
	\begin{theorem}\label{consta from negacyclic}
		Let $\mathtt{C_{_{\theta}}}$ be a $\beta_{_{\theta}}$- constacyclic code over $\mathtt{R_{_{\theta}}}, \theta \in \{0,1,\nu\}$ of arbitrary length $N=n.2^{s},$ where $s$ and $n$ are positive integers and $n$ is odd. Then $\mathtt{C_{_{\theta}}}= \langle l_{_{\theta_{1}}}(z),l_{_{\theta_{2}}}(z)\rangle,$ where $l_{_{\theta_{i}}}(z) \in \dfrac{\mathtt{R_{_{\theta}}}[z]}{\langle z^{N}-\beta_{_{\theta}} \rangle}$ for $ i=1,2.$ Further, $l_{_{\theta_{1}}}(z)=  (z^{n}-1)^{t_{1}}+k_{_{\theta}}h(z)$ and $ l_{_{\theta_{2}}}(z)=k_{_{\theta}}(z^{n}-1)^{t_{2}},$ where $0 \leq t_{1},t_{2} \leq 2^{s+1}$ and $h(z)$ is in $\dfrac{Z_{4}[z]}{\langle z^{N}+1 \rangle}.$
	\end{theorem}
	\begin{proof}
		Let $\mathtt{C}_{_{\theta}}$ be a $\beta_{_{\theta}}$- constacyclic code over $\mathtt{R}_{_{\theta}}$, $\theta \in \{0,1,\nu\}$ of length $N=n.2^{s},$ where $s$ and $n$ are positive integers and $n$ is odd. Let $\phi_{_{\theta}}$ be the ring homomorphism as defined in Lemma \ref{iso homo lemma} for $\theta \in \{0,1,\nu\}.$ It is easy to see that, $\phi_{_{\theta}}(\mathtt{C}_{_{\theta}})$ is a negacyclic code of length $N$ over $Z_{4}.$ Using Lemma \ref{nega cyclic gen over $Z_{4}$}, we get $\phi_{_{\theta}}(\mathtt{C}_{_{\theta}})$ = $\langle (z^{n}-1)^{t_{1}}\rangle$, where $0 \leq t_{1} \leq 2^{s+1}.$ 
		Also $ker_{_{\theta}}=\langle k_{_{\theta}} \rangle.$ Consider $I=\biggl\{a(z)\in \dfrac{Z_{4}[z]}{\langle z^{N}+1\rangle}: k_{_{\theta}}a(z)\in ker_{_{\theta}}\biggr\}.$ Clearly, $I$ is an ideal of $\dfrac{Z_{4}[z]}{\langle z^{N}+1\rangle}.$ Thus, $ker_{_{\theta}}$ is an ideal of $k_{_{\theta}}Z_{4}[z]/\langle z^{
			N}+1\rangle.$ Therefore, $ker_{_{\theta}}$ is of the form $k_{_{\theta}}J,$ where $J$ is a negacyclic code of length $N$ over $Z_{4}.$ Therefore, using Lemma \ref{nega cyclic gen over $Z_{4}$}, we obtain $ker_{_{\theta}}$ = $k_{_{\theta}} \langle (z^{n}-1)^{t_{2}}\rangle$, where $0 \leq t_{2} \leq 2^{s+1}.$
		
		It follows that
		$\mathtt{C}_{_{\theta}}$ is generated by $l_{_{\theta_{1}}}(z)= (z^{n}-1)^{t_{1}}+k_{_{\theta}}h(z)$ and $l_{_{\theta_{2}}}(z)=k_{_{\theta}}(z^{n}-1)^{t_{2}},$ where $0 \leq t_{1},t_{2} \leq 2^{s+1}$ such that $l_{_{\theta_{i}}}(z) \in \dfrac{\mathtt{R_{_{\theta}}}[z]}{\langle z^{N}-\beta_{_{\theta}} \rangle}$ for $i=1,2$ and $h(z)$ is in $\dfrac{Z_{4}[z]}{\langle z^{N}+1\rangle}.$  
	\end{proof}
	
	\begin{remark}
		It can be easily seen that for a $\beta_{_{\theta}}$- constacyclic code $\mathtt{C_{_{\theta}}}$ of arbitrary length $N=n.2^{s},$ (where $s$ and $n$ are positive integers and $n$ is odd) over $\mathtt{R_{_{\theta}}}, \theta \in \{0,1,\nu\},$ the negacyclic codes Res$(\mathtt{C}_{_{\theta}})$ and Tor$(\mathtt{C}_{_{\theta}})$ over $Z_{4}$ are generated respectively by the polynomials $(z^{n}-1)^{t_{1}}$ and $(z^{n}-1)^{t_{2}}.$ Thus  $(z^{n}-1)^{t_{1}}$ and $(z^{n}-1)^{t_{2}},$ where $0 \leq t_{1},t_{2} \leq 2^{s+1}$ can be considered to be the unique minimal degree polynomial generators of the principally generated codes Res$(\mathtt{C}_{_{\theta}})$ and Tor$(\mathtt{C}_{_{\theta}})$ respectively.   
	\end{remark}
	
	\begin{theorem} \label{theorem str 2 cc beta}	Let $\mathtt{C}_{_{\theta}}$ be a 
		$\beta_{_{\theta}}$- constacyclic code of arbitrary length $N=n.2^{s}$ over the ring $\mathtt{R}_{_{\theta}}, \theta \in \{0,1,\nu\}.$
		 Then there exists a set of generators  $\{t_{_{\theta_{1}}}(z),t_{_{\theta_{2}}}(z)\}$ of $\mathtt{C}_{_{_\theta}},$ where $t_{_{\theta_{1}}}(z)=(z^{n}-1)^{t_{1}}+k_{_{\theta}}t(z)$ and $t_{_{\theta_{2}}}(z)=k_{_{\theta}}(z^{n}-1)^{t_{2}},$ where $0 \leq t_{1},t_{2} \leq 2^{s+1}$ such that $t_{_{\theta_{i}}}(z) \in \dfrac{\mathtt{R_{_{\theta}}}[z]}{\langle z^{N}-\beta_{_{\theta}} \rangle}$ for $ i=1,2$ and $t(z)$ is in $\dfrac{Z_{4}[z]}{\langle z^{N}+1\rangle}.$ Further, $(z^{n}-1)^{t_{1}}$ and $(z^{n}-1)^{t_{2}}$ are unique minimal degree polynomial generators of Res$(\mathtt{C}_{_{\theta}})$ and Tor$(\mathtt{C}_{_{\theta}})$ respectively. 
		Also, either $t(z)=0$ or deg $t(z)<$  deg $(z^{n}-1)^{t_{2}}=nt_{2}.$ 
	\end{theorem}
	\begin{proof} 
			Let $\mathtt{C}_{_{\theta}}$ be a 
		$\beta_{_{\theta}}$- constacyclic code of arbitrary length $N=n.2^{s}$ over the ring $\mathtt{R}_{_{\theta}}, \theta \in \{0,1,\nu\}.$ The by Theorem \ref{consta from negacyclic} above $\mathtt{C}_{_{\theta}}
		= \langle l_{_{\theta_{1}}}(z),l_{_{\theta_{2}}}(z)\rangle,$ where 
		$l_{_{\theta_{1}}}(z)=(z^{n}-1)^{t_{1}}+k_{_{\theta}}h(z)$ and $l_{_{\theta_{2}}}(z)=k_{_{\theta}}(z^{n}-1)^{t_{2}}.$ If $h(z)=0$ or deg $h(z) <$ $nt_{2},$ then we get the required result. Otherwise, let us suppose that deg $h(z) \geq$ deg $(z^{n}-1)^{t_{2}}.$ By division algorithm, there exist some polynomials $q(z), t(z) \in Z_{4}[z]$ such that $h(z)=q(z)(z^{n}-1)^{t_{2}}+t(z),$ where either $t(z)=0$ or deg $t(z) <$ $nt_{2}.$ Now consider, the polynomial  $t_{_{\theta_{1}}}(z)= l_{_{\theta_{1}}}(z)-q(z)l_{_{\theta_{2}}}(z)=(z^{n}-1)^{t_{1}}+k_{_{\theta}}t(z),$ where either $t(z)=0$ or deg $t(z) <$ $nt_{2}.$ Clearly, $t_{_{\theta_{1}}}(z) \in \mathtt{C_{_{\theta}}}.$
		Since $t_{_{\theta_{1}}}(z)$ is a linear combination of $l_{_{\theta_{1}}}(z)$ and $l_{_{\theta_{2}}}(z),$ we have $\mathtt{C}_{_{\theta}} = \langle l_{_{\theta_{1}}}(z),l_{_{\theta_{2}}}(z)\rangle = \langle t_{_{\theta_{1}}}(z),t_{_{\theta_{2}}}(z)\rangle.$  This completes the proof of the theorem. 
	\end{proof}
	In the following theorem, it is proved that the set of generators $\{t_{_{\theta_{1}}}(z),t_{_{\theta_{2}}}(z)\}$ of a $\beta_{_{\theta}}$- constacyclic code $\mathtt{C_{_{\theta}}}$ of arbitrary length over $\mathtt{R}_{_{\theta}}, \theta \in \{0,1,\nu\}$ obtained in Theorem \ref{theorem str 2 cc beta} are uniquely determined. 
	\begin{theorem} \label{theorem str 3 cc beta}	Let $\mathtt{C}_{_{\theta}}=\langle t_{_{\theta_{1}}}(z),t_{_{\theta_{2}}}(z)\rangle$ be a $\beta_{_{\theta}}$- constacyclic code of arbitrary length $N=n.2^{s},$ (where $s$ and $n$ are positive integers and $n$ is odd) over the ring $\mathtt{R}_{_{\theta}}, \theta \in \{0,1,\nu\}$,  where  $t_{_{\theta_{i}}}(z) \in \dfrac{\mathtt{R_{_{\theta}}}[z]}{\langle z^{N}-\beta_{_{\theta}} \rangle}$ for $1 \leq i \leq 2$ are polynomials as obtained in Theorem \ref{theorem str 2 cc beta}.
		Then the polynomials $t_{_{\theta_{i}}}(z)$ for $ 1 \leq i \leq 2$ are uniquely determined.
	\end{theorem}
	
	\begin{proof} If possible, consider another set of generators  $\{r_{_{\theta_{1}}}(z),r_{_{\theta_{2}}}(z)\}$ of $\mathtt{C}_{_{\theta}}$, where $r_{_{\theta_{1}}}(z)=  (z^{n}-1)^{t_{1}}+k_{_{\theta}}r(z)$ and $ r_{_{\theta_{2}}}(z)=k_{_{\theta}}(z^{n}-1)^{t_{2}},$ where $0 \leq t_{1},t_{2} \leq 2^{s+1}$ and $r(z)$ is in $\dfrac{Z_{4}[z]}{\langle z^{N}+1 \rangle}$ with either $r(z)=0$ or deg $r(z) <$ $nt_{2}.$ 
		Also, $(z^{n}-1)^{t_{1}}$ and $(z^{n}-1)^{t_{2}}$ are unique minimal degree polynomial generators of Res$(\mathtt{C}_{_{\theta}})$ and Tor$(\mathtt{C}_{_{\theta}})$ respectively.
		
		We have, $t_{_{\theta_{1}}}(z)-r_{_{\theta_{1}}}(z)=k_{_{\theta}}(t(z)-r(z))\in \mathtt{C}_{_{\theta}}.$
		This implies that $(t(z)-r(z)) \in$ Tor$(\mathtt{C}_{_{\theta}}) =\langle (z^{n}-1)^{t_{2}}\rangle.$ Also, deg $(t(z)-r(z))<$ deg $(z^{n}-1)^{t_{2}}$ which is a contradiction because $(z^{n}-1)^{t_{2}}$ is a minimal degree polynomial in  Tor($\mathtt{C}_{_{\theta}})$. It follows that $t(z)=r(z).$ Hence, the uniqueness of the polynomials $t_{_{\theta_{i}}}(z)$ for $ 1 \leq i \leq 2$ is established.
	\end{proof}	
	
	\section{Rank and cardinality of a constacyclic code of arbitrary length over $\mathtt{R}_{_{\theta}}$}
	In this section, the rank and cardinality of a constacyclic code of arbitrary length over $\mathtt{R}_{_{\theta}}, \theta \in \{0,1,\nu\}$ have been obtained by determining a minimal spanning set of the code.
	
	The rank and cardinality of a cyclic code of arbitrary length over a non-chain ring of the type $Z_{4}+\nu Z_{4}$ where $\nu^{2} \in Z_{4}+\nu Z_{4}$ have been obtained by Jain et al. \cite{self paper}. It is observed that a similar approach can be used to obtain the rank and cardinality of an $\alpha_{_{\theta}}$- constacyclic code of arbitrary length over a non-chain ring of the type $Z_{4}+\nu Z_{4}$ where $\nu^{2} \in Z_{4}+\nu Z_{4}.$
	
	In the following theorem, the rank and cardinality of an $\alpha_{_{\theta}}$- constacyclic code of arbitrary length over $\mathtt{R}_{_{\theta}}, \theta \in \{0,1,\nu\}$ have been obtained by determining a minimal spanning set of the code.
	\begin{theorem} \label{theorem minimal spanning cc} Let $\mathtt{C}_{_{\theta}}=\langle t_{_{\theta_{1}}}(z),t_{_{\theta_{2}}}(z),t_{_{\theta_{3}}}(z),t_{_{\theta_{4}}}(z)\rangle$ be an $\alpha_{_{\theta}}$- constacyclic code of arbitrary length $n$ over the ring $\mathtt{R}_{_{\theta}},\theta \in \{0,1,\nu\}$, where $t_{_{\theta_{i}}}(z) \in \dfrac{\mathtt{R_{_{\theta}}}[z]}{\langle z^{n}-\alpha_{_{\theta}} \rangle}$ for $1 \leq i \leq 4$. Further, let $t_{_{\theta_{1}}}(z)= t_{_{11}}(z)+2t_{_{12}}(z)+k_{_{\theta}}t_{_{13}}(z)+2k_{_{\theta}}t_{_{14}}(z)$, $t_{_{\theta_{2}}}(z)=2t_{_{22}}(z)+k_{_{\theta}}t_{_{23}}(z)+2k_{_{\theta}}t_{_{24}}(z)$, $t_{_{\theta_{3}}}(z)=k_{_{\theta}}t_{_{33}}(z)+2k_{_{\theta}}t_{_{34}}(z)$ and $t_{_{\theta_{4}}}(z)=2k_{_{\theta}}t_{_{44}}(z)$ be in the unique form as given in Theorem \ref{theorem str 3cc}. Let $s_{_{i}}=$ deg $t_{_{ii}}(z)$ for $1 \leq i \leq 4$ and $\tilde{s}=	min\{s_{_{2}},s_{_{3}}\}.$  Then $$rank(\mathtt{C}_{_{\theta}})= n+s_{_{1}}+\tilde{s}-s_{_{2}}-s_{_{3}}-s_{_{4}}$$ and   \[ \vert\mathtt{C_{_{\theta}}}\vert= \begin{cases}
			2^{4n+s_{_{1}}+\tilde{s}-3s_{_{2}}-2s_{_{3}}-s_{_{4}}}& ;~t_{_{23}}(z) \ne 0 \\
			2^{4n+\tilde{s}-2s_{_{2}}-2s_{_{3}}-s_{_{4}}}& ;~t_{_{23}}(z) =0 \\
		\end{cases}. \]
	\end{theorem}
	\begin{proof}
		The proof follows on the same lines as the proof of Theorem 11 \cite{self paper}.
	\end{proof}
	
			
			In the following theorem, the rank and cardinality of a $\beta_{_{\theta}}$- constacyclic code of arbitrary length over the ring $\mathtt{R}_{_{\theta}}, \theta \in \{0,1,\nu\}$ have been obtained by determining a minimal spanning set of the code. 
			
			\begin{theorem} \label{theorem minimal spanning cc beta} Let $\mathtt{C}_{_{\theta}}=\langle t_{_{\theta_{1}}}(z),t_{_{\theta_{2}}}(z)\rangle$ be a  $\beta_{_{\theta}}$- constacyclic code of arbitrary length $N=n.2^{s},$ (where $s$ and $n$ are positive integers and $n$ is odd) over the ring $\mathtt{R}_{_{\theta}},\theta \in \{0,1,\nu\}.$ Further, let the generators $t_{_{\theta_{1}}}(z)= (z^{n}-1)^{t_{1}}+k_{_{\theta}}t(z)$ and $ t_{_{\theta_{2}}}(z)=k_{_{\theta}}(z^{n}-1)^{t_{2}},$ where $0 \leq t_{1},t_{2} \leq 2^{s+1}$ and $t(z) \in \dfrac{Z_{4}[z]}{\langle z^{N}+1 \rangle}$ be in the unique form as given in Theorem \ref{theorem str 3 cc beta}. Then $rank(\mathtt{C}_{_{\theta}}) = N-nt_{2}$ and $\vert \mathtt{C}_{_{\theta}}\vert =
				2^{4N-2nt_{1}-2nt_{2}}.$  
				
			\end{theorem}
			\begin{proof} It can be easily seen that the set $\mathtt{A}_{_{\theta}}=\{t_{_{\theta_{1}}}(z),zt_{_{\theta_{1}}}(z),\cdots,z^{N-nt_{1}-1}t_{_{\theta_{1}}}(z),\\t_{_{\theta_{2}}}(z),zt_{_{\theta_{2}}}(z),\cdots,z^{N-nt_{2}-1}t_{_{\theta_{2}}}(z)\}$ is a spanning set of $\mathtt{C}_{_{\theta}}.$
				
				To prove that $rank~(\mathtt{C}_{_{\theta}})$ is $N-nt_{2},$ it is sufficient to show that the set $\mathtt{B}_{_{\theta}}=\{ t_{_{\theta_{1}}}(z),zt_{_{\theta_{1}}}(z),\cdots,z^{N-nt_{1}-1}t_{_{{\theta}_{1}}}(z),t_{_{\theta_{2}}}(z),zt_{_{\theta_{2}}}(z),\cdots,z^{nt_{1}-nt_{2}-1}t_{_{\theta_{2}}}(z)\}$ is a minimal spanning set of $\mathtt{C}_{_{\theta}}.$
				
				To prove that the set $\mathtt{B}_{_{\theta}}$ spans $\mathtt{C}_{_{\theta}},$ it is enough to show that $z^{nt_{1}-nt_{2}}t_{_{\theta_{2}}}(z) \in span({\mathtt{B}_{_{\theta}}})$. 
				Since deg $z^{nt_{1}-nt_{2}} t_{_{\theta_{2}}}(z)$ = deg $t_{_{\theta_{1}}}(z)=nt_{1},$ there exists a polynomial $r_{_{1}}(z)$ such that 
				\begin{equation}\label{eqn8 cc beta}
					r_{_{1}}(z)=z^{nt_{1}-nt_{2}}t_{_{\theta_{2}}}(z)-k_{_{\theta}}t_{_{\theta_{1}}}(z).
				\end{equation}
				Clearly, $r_{_{1}}(z) \in \mathtt{C_{_{\theta}}}.$ Moreover, either $r_{_{1}}(z)=0$ or  deg $r_{_{1}}(z)< nt_{1}.$ If $r_{_{1}}(z)=0,$ then  $z^{nt_{1}-nt_{2}}t_{_{\theta_{2}}}(z) \in span(\mathtt{B_{_{\theta}}}).$ If deg $
				r_{_{1}}(z) < nt_{1},$ then it is easy to see that $r_{_{1}}(z)$ is of the type $t_{_{\theta_{2}}}(z).$ Due to the minimality of degree of $t_{_{\theta_{2}}}(z),$ we have deg $r_{_{1}}(z) \geq nt_{2}.$ Therefore, there exists a polynomial $r_{_{2}}(z)$ such that $r_{_{2}}(z)=r_{_{1}}(z)-z^{{\text{ deg }r_{_{1}}(z)-nt_{2}}}t_{_{\theta_{2}}}(z).$ It is easy to see that $r_{_{2}}(z) \in \mathtt{C_{_{\theta}}}$ and it is of the type $t_{_{\theta_{2}}}(z).$ Also, either $r_{_{2}}(z)=0$ or deg $ r_{_{2}}(z) <$ deg $r_{_{1}}(z).$ If $r_{_{2}}(z)=0,$ then $r_{_{1}}(z)=z^{{\text{ deg }r_{_{1}}(z)-nt_{2}}}t_{_{\theta_{2}}}(z).$ Subsituting the value of $r_{_{1}}(z)$ in equation (\ref{eqn8 cc beta}), we see that $z^{nt_{1}-nt_{2}} t_{_{\theta_{2}}}(z) \in span(\mathtt{B_{_{\theta}}}).$ If deg $ r_{_{2}}(z) <$ deg $r_{_{1}}(z),$ then after repeating the argument a finite number of times we obtain a polynomial $r_{_{l}}(z) = r_{_{l-1}}(z)- z^{{\text{ deg }r_{_{l-1}}(z)-nt_{2}}}t_{_{\theta_{2}}}(z)$ such that $r_{_{l}}(z) \in \mathtt{C_{_{\theta}}}$ and it is of the type $t_{_{\theta_{2}}}(z).$ Moreover, $r_{_{l}}(z)=0$ or deg $r_{_{l}}(z) < nt_{2}.$ Since $r_{_{l}}(z)$ is of the type $t_{_{\theta_{2}}}(z),$ deg $r_{_{l}}(z)$ cannot be less than $nt_{2}.$ Therefore, $r_{_{l}}(z)=0.$ Hence, from equation (\ref{eqn8 cc beta}), we have,
				$z^{nt_{1}-nt_{2}} t_{_{\theta_{2}}}(z)=k_{_{\theta}}t_{_{\theta_{1}}}(z)+r_{_{1}}(z)=k_{_{\theta}}t_{_{\theta_{1}}}(z)+z^{{\text{ deg }r_{_{1}}(z)-nt_{2}}}t_{_{\theta_{2}}}(z)+r_{_{2}}(z)=k_{_{\theta}}t_{_{\theta_{1}}}(z)+z^{{\text{ deg }r_{_{1}}(z)-nt_{2}}}g_{_{\theta_{2}}}(z)+z^{{\text{ deg }r_{_{2}}(z)-nt_{2}}}t_{_{\theta_{2}}}(z)+\cdots+z^{{\text{ deg }r_{_{l-1}}(z)-nt_{2}}}t_{_{\theta_{2}}}(z).$ It follows that $z^{nt_{1}-nt_{2}} t_{_{\theta_{2}}}(z) \in span(\mathtt{B_{_{\theta}}}).$ 
				
				Now to prove that the set $\mathtt{B_{_{\theta}}}$ is a minimal spanning set, it is enough to show that none of $z^{N-nt_{1}-1}t_{_{\theta_{1}}}(z)$ and $ z^{nt_{1}-nt_{2}-1}t_{_{\theta_{2}}}(z)$ can be written as a linear combination of other elements of $\mathtt{B}_{_{\theta}}.$ Suppose, if possible, that $z^{N-nt_{1}-1}t_{_{\theta_{1}}}(z)$ can be written as a linear combination of other elements of $\mathtt{B}_{_{\theta}},$ i.e, 
				\begin{equation}\label{equation 2.9 cc beta}
					z^{N-nt_{1}-1}t_{_{\theta_{1}}}(z)= a(z)t_{_{\theta_{1}}}(z)
					+b(z)t_{_{\theta_{2}}}(z),
				\end{equation}
				where deg $a(z) < N-nt_{1}-1$ and deg $b(z) < nt_{1}-nt_{2}.$  Multiplying equation (\ref{equation 2.9 cc beta}) on both sides by $2k_{_{\theta}}$ for $\theta \in \{0,1\},$ we get
				\begin{equation}\label{equation 2.10 cc beta}
					2k_{_{\theta}}z^{N-nt_{1}-1} (z^{n}-1)^{t_{1}}= 2k_{_{\theta}} a(z)(z^{n}-1)^{t_{1}}, ~\theta \in \{0,1\}.
				\end{equation}
				Multiplying equation (\ref{equation 2.9 cc beta}) on both sides by $2(k_{_{\theta}}-1)$ for $\theta =\nu,$ we get
				\begin{equation}\label{equation 2.11 cc beta}
					2(k_{_{\theta}}-1)z^{N-nt_{1}-1} (z^{n}-1)^{t_{1}}= 2(k_{_{\theta}}-1)a(z)(z^{n}-1)^{t_{1}}, ~\theta =\nu.
				\end{equation}
				The equations (\ref{equation 2.10 cc beta}) and (\ref{equation 2.11 cc beta}) are not possible as degrees of left hand side and right hand side in each of these equations do not match. Thus, $ z^{N-nt_{1}-1}t_{_{\theta_{1}}}(z)$ can not be written as a linear combination of other elements of $\mathtt{B_{_{\theta}}}.$ Using a similar argument, it can be shown that $z^{nt_{1}-nt_{2}-1}t_{_{\theta_{2}}}(z)$ can not be written as a linear combination of other elements of $\mathtt{B}_{_{\theta}}.$ Hence, $\mathtt{B}_{_{\theta}}$ is a minimal spanning set of $\mathtt{C}_{_{\theta}}.$ 
				Further,  $rank(\mathtt{C}_{_{\theta}})=$ Number of elements in $\mathtt{B}_{_{\theta}}=(N-nt_{1})+(nt_{1}-nt_{2})=N-nt_{2}.$  It can be easily seen that cardinality of $\mathtt{C}_{_{\theta}}$ is  equal to $
				2^{4N-2nt_{1}-2nt_{2}}.$  
			\end{proof}
	
	\section{Reversible constacyclic codes over $\mathtt{R_{_{\theta}}}$}
	In this section, necessary and sufficient conditions for a constacyclic code of arbitrary length over $\mathtt{R}_{_{\theta}}, \theta \in \{0,1,\nu\}$ to be reversible have been obtained. We begin with some definitions and results that we require to proceed further in this section.
	
	Consider a finite commutative ring $\mathtt{R}$ with unity and a linear code $\mathtt{C}$ over $\mathtt{R}.$ A code $\mathtt{C}$ of length $n$ is said to be a reversible code if for every $\mathtt{s}=(\mathtt{s}_{_{0}},\mathtt{s}_{_{1}},\cdots,\mathtt{s}_{{n-1}})$ in $\mathtt{C}$, the codeword $\mathtt{s^{r}}= (\mathtt{s}_{_{n-1}},\mathtt{s}_{_{n-2}},\cdots,\mathtt{s}_{_{0}})$ also belongs to $\mathtt{C}$. The polynomial representation of the codeword $\mathtt{s^{r}}$ is $\mathtt{s}_{_{n-1}}+\mathtt{s}_{_{n-2}}z+\cdots+\mathtt{s}_{_{0}}z^{n-1}$ which is denoted by $\mathtt{s^{r}}(z).$ For a polynomial $g(z)$ of degree $t \leq n-1$, $g^{*}(z)= z^{t}g(z^{-1})$ is defined as its reciprocal polynomial and $g^{\mathtt{r}}(z)=z^{n-t-1}g^{*}(z)$ is defined as its reverse polynomial. A polynomial $g(z)$ is said to be self reciprocal if and only if $g^{*}(z)=g(z)$. 
	
	\begin{lemma} {\label{lemma poly add and multi for rev}} (\cite{jasbir thesis}) Let $g_{_{1}}(z)$ and $g_{_{2}}(z)$ be any two polynomials in $\mathtt{R}[z]$ with deg $g_{_{1}}(z)\geq$ deg $g_{_{2}}(z)$. Then
		\begin{enumerate}
			\item [$(i)$]$\big(g_{_{1}}(z)+g_{_{2}}(z)\big)^{*}= g_{_{1}}^{*}(z)+z^{i}g_{_{2}}^{*}(z)$, where $i= $deg $g_{_{1}}(z)-$ deg $g_{_{2}}(z),$
			\item [$(ii)$] $\big(g_{_{1}}(z)g_{_{2}}(z)\big)^{*}= g_{_{1}}^{*}(z)g_{_{2}}^{*}(z).$
		\end{enumerate}
	\end{lemma}
		It has been proved by J.Kaur \cite{jasbir thesis} that for a cyclic code $\mathtt{C}$ over $\mathtt{R}$ with generators $g_{_{1}}(z),g_{_{2}}(z),\cdots,g_{_{k}}(z)$ is a reversible cyclic code if and only if $g_{_{i}}^{*}(z) \in \mathtt{C}$ for all $ 1\leq i\leq k$. We observe that this result can be extended to constacylic codes over a finite commutative ring $\mathtt{R}$ with unity. We state and prove below the corresponding result which gives necessary and sufficient conditions for a $\lambda$- constacyclic code over $\mathtt{R}$ to be reversible in terms of its generators. 
		
		\begin{lemma}\label{jasbir result on consta}
			Let $\mathtt{C}$ be a $\lambda$- constacyclic code over $\mathtt{R}$ with generators $g_{_{1}}(z),g_{_{2}}(z),\\\cdots,g_{_{k}}(z)$. Then $\mathtt{C}$ is a reversible $\lambda$- constacyclic code if and only if $g_{_{i}}^{*}(z) \in \mathtt{C}$ for all $ 1\leq i\leq k$.
		\end{lemma}
		\begin{proof}
			Let $\mathtt{C}$ be a $\lambda$- constacyclic code over $\mathtt{R}.$ Then $\mathtt{C}$ can be viewed as an ideal of $\dfrac{\mathtt{R}[z]}{\langle z^{n}-\lambda \rangle}.$ Let $\mathtt{C}= \langle g_{_{1}}(z),g_{_{2}}(z),\cdots,g_{_{k}}(z) \rangle.$ Suppose $g_{_{i}}^{*}(z) \in \mathtt{C}$ for all $ 1\leq i\leq k.$ Let $\mathtt{s}(z)$ be an arbitrary polynomial in $\mathtt{C}.$ Then 
			\begin{equation*}
				\mathtt{s}(z)=g_{_{1}}(z)h_{_{1}}(z)+g_{_{2}}(z)h_{_{2}}(z)+\cdots+g_{_{k}}(z)h_{_{k}}(z),
			\end{equation*}
			where $h_{_{i}}(z) \in \mathtt{R}[z]$ for $ 1\leq i\leq k.$ This implies 
			\begin{equation*}
				\mathtt{s}^{*}(z)=z^{\lambda_{_{1}}}g_{_{1}}^{*}(z)h_{_{1}}^{*}(z)+z^{\lambda_{_{2}}}g_{_{2}}^{*}(z)h_{_{2}}^{*}(z)+\cdots+z^{\lambda_{_{k}}}g_{_{k}}^{*}(z)h_{_{k}}^{*}(z),
			\end{equation*}
			where $\lambda_{_{i}} \in \mathbb{Z}^{+} \cup \{0\}$ and $h_{_{i}}(z) \in \mathtt{R}[z]$ for $ 1\leq i\leq k.$ As $g_{_{i}}^{*}(z) \in \mathtt{C}$ for $ 1\leq i\leq k,$ it follows that $\mathtt{s}^{*}(z) \in \mathtt{C}$ because $\mathtt{C}$ is an ideal of $\dfrac{\mathtt{R}[z]}{\langle z^{n}-\lambda \rangle}.$ Therefore, $\mathtt{s}^{\mathtt{r}}(z)=z^{n-t-1}\mathtt{s}^{*}(z) \in \mathtt{C},$ where $t$ is the degree of $\mathtt{s}(z).$ Thus $\mathtt{C}$ is a reversible $\lambda$- constacyclic code.	
			
			Conversely, let $\mathtt{C}$ be a reversible $\lambda$- constacyclic code. Then $g_{_{i}}^{\mathtt{r}}(z) \in \mathtt{C}$ for all $1\leq i\leq k.$ Therefore, $g_{_{i}}^{*}(z)=z^{t-n+1}g_{_{i}}^{\mathtt{r}}(z) \in \mathtt{C}$ for all $ 1\leq i\leq k.$ This proves the result.
		\end{proof}
		
		In the following theorem, it is shown that the torsion code of a reversible $u_{_{\theta}}$- constayclic code $\mathtt{C}_{_{\theta}}$ of arbitrary length over $\mathtt{R}_{_{\theta}}, \theta \in \{0,1,\nu\}$ is a reversible code over $Z_{4}.$
		
		\begin{theorem}\label{tor consta from cyclic}
			Let $\mathtt{C_{_{\theta}}}$ be a $u_{_{\theta}}$- constacyclic code of arbitrary length $n$ over $\mathtt{R_{_{\theta}}}, \theta \in \{0,1,\nu\}.$ Then
			\begin{enumerate}
				\item [$(a)$]  Tor$(\mathtt{C}_{_{\theta}})$ is a reversible cyclic code over $Z_{4}$ if $u_{_{\theta}}$ is of the type $\alpha_{_{\theta}}.$
				\item [$(b)$]  Tor$(\mathtt{C}_{_{\theta}})$ is a reversible negacyclic code over $Z_{4}$ if $u_{_{\theta}}$ is of the type $\beta_{_{\theta}}.$
			\end{enumerate}
		\end{theorem}
		\begin{proof} $(a)$ Let $\mathtt{C}_{_{\theta}}$ be a reversible $\alpha_{_{\theta}}$- constacyclic code of arbitrary length $n$ over $\mathtt{R}_{_{\theta}}, \theta \in \{0,1,\nu\}.$ Then $\mathtt{C}_{_{\theta}}=\langle h_{_{\theta_{1}}}(z),h_{_{\theta_{2}}}(z),h_{_{\theta_{3}}}(z),h_{_{\theta_{4}}}(z)\rangle,$ where $h_{_{\theta_{i}}}(z) \in \dfrac{\mathtt{R_{_{\theta}}}[z]}{\langle z^{n}-\alpha_{_{\theta}} \rangle}$ for $1 \leq i \leq 4$. Further, $h_{_{\theta_{1}}}(z)= h_{_{11}}(z)+2h_{_{12}}(z)+k_{_{\theta}}h_{_{13}}(z)+2k_{_{\theta}}h_{_{14}}(z)$, $h_{_{\theta_{2}}}(z)=2h_{_{22}}(z)+k_{_{\theta}}h_{_{23}}(z)+2k_{_{\theta}}h_{_{24}}(z)$, $h_{_{\theta_{3}}}(z)=k_{_{\theta}}h_{_{33}}(z)+2k_{_{\theta}}h_{_{34}}(z)$ and $h_{_{\theta_{4}}}(z)=2k_{_{\theta}}h_{_{44}}(z)$ such that the polynomials $h_{_{ij}}(z)$ are in $Z_{2}[z]/{\langle z^{n}-1\rangle}$ for $1 \leq i \leq j \leq 4.$ It is easy to see that Tor$(\mathtt{C}_{_{\theta}})$ is a cyclic code generated by $h_{_{33}}(z)+2h_{_{34}}(z)$ and $2h_{_{44}}(z)$ over $Z_{4}. $
			Since $k_{_{\theta}}(h_{_{33}}(z)+2h_{_{34}}(z)) \in \mathtt{C}_{_{\theta}}$ and $\mathtt{C}_{_{\theta}}$ is reversible $\alpha_{_{\theta}}$- constacyclic code, therefore, $k_{_{\theta}}(h_{_{33}}(z)+2h_{_{34}}(z))^{*} \in \mathtt{C}_{_{\theta}}$ by Lemma \ref{jasbir result on consta}. It follows that $(h_{_{33}}(z)+2gh{_{34}}(z))^{*} \in$ Tor$(\mathtt{C}_{_{\theta}})$. Similarly, $(2h_{_{44}}(z))^{*} \in$ Tor$(\mathtt{C}_{_{\theta}}).$ Hence, Tor$(\mathtt{C}_{_{\theta}})=\langle h_{_{33}}(z)+2h_{_{34}}(z),2h_{_{44}}(z)\rangle$ is a reversible cyclic code over $Z_{4}$ by Lemma \ref{jasbir result on consta}.
			
			$(b)$ Let $\mathtt{C_{_{\theta}}}$ be a reversible $\beta_{_{\theta}}$- constacyclic code of length $N=n.2^{s}$ over $\mathtt{R}_{_{\theta}}, \theta \in \{0,1,\nu\}.$ Then $\mathtt{C_{_{\theta}}}= \langle l_{_{\theta_{1}}}(z),l_{_{\theta_{2}}}(z)\rangle,$ where $l_{_{\theta_{i}}}(z) \in \dfrac{\mathtt{R_{_{\theta}}}[z]}{\langle z^{N}-\beta_{_{\theta}} \rangle}$ for $ i=1,2.$ Further, $l_{_{\theta_{1}}}(z)=  (z^{n}-1)^{t_{1}}+k_{_{\theta}}h(z)$ and $ l_{_{\theta_{2}}}(z)=k_{_{\theta}}(z^{n}-1)^{t_{2}}$ where $0 \leq t_{1},t_{2} \leq 2^{s+1}$ and $h(z)$ is in $\dfrac{Z_{2}[z]}{\langle z^{N}+1 \rangle}.$ It is easy to see that Tor$(\mathtt{C}_{_{\theta}})$ is a negacyclic code generated by $(z^{n}-1)^{t_{2}}$ over $Z_{4}.$ Since $k_{_{\theta}}(z^{n}-1)^{t_{2}} \in \mathtt{C}_{_{\theta}}$ and $\mathtt{C}_{_{\theta}}$ is a reversible $\beta_{_{\theta}}$- constacyclic code, therefore, $k_{_{\theta}}\big((z^{n}-1)^{t_{2}}\big)^{*} \in \mathtt{C}_{_{\theta}}$ by Lemma \ref{jasbir result on consta}. It follows that $\big((z^{n}-1)^{t_{2}}\big)^{*} \in$ Tor$(\mathtt{C}_{_{\theta}}).$ Hence, Tor$(\mathtt{C}_{_{\theta}})=\langle (z^{n}-1)^{t_{2}} \rangle$ is a reversible negacyclic code over $Z_{4}$ by Lemma \ref{jasbir result on consta}. 
		\end{proof}
		
		Necessary and sufficient conditions for a cyclic code of arbitrary length over a non-chain ring of the type $Z_{4}+\nu Z_{4},$ where $\nu^{2} \in Z_{4}+\nu Z_{4}$ to be reversible have been obtained by Jain et al. \cite{arxiv rev and rev comp}. It is observed that a similar approach can be used to obtain necessary and sufficient conditions of an $\alpha_{_{\theta}}$- constacyclic code of arbitrary length over a non-chain ring of the type $Z_{4}+\nu Z_{4},$ where $\nu^{2} \in Z_{4}+\nu Z_{4}$ to be reversible.
		
		The following theorem gives necessary and sufficient conditions for an $\alpha_{_{\theta}}$- constacyclic code $\mathtt{C}_{_{\theta}}$ of arbitrary length $n$ over $\mathtt{R}_{_{\theta}}, \theta \in \{0,1,\nu\}$ to be a reversible code.
		
		\begin{theorem}
			{\label{const from cyc theorem NASC rev}} Let $\mathtt{C}_{_{\theta}}$ be an $\alpha_{_{\theta}}$- constacyclic code of arbitrary length $n$ over $\mathtt{R}_{_{\theta}}, \theta \in \{0,1,\nu\}$ generated by $h_{_{\theta_{1}}}(z),h_{_{\theta_{2}}}(z),h_{_{\theta_{3}}}(z)$ and $h_{_{\theta_{4}}}(z),$ where $h_{_{\theta_{i}}}(z) \in \dfrac{\mathtt{R_{_{\theta}}}[z]}{\langle z^{n}-\alpha_{_{\theta}} \rangle}$ for $1 \leq i \leq 4$ are as obtained in Theorem \ref{theorem str 3cc}. Then $\mathtt{C}_{_{\theta}}$ is a
			reversible $\alpha_{_{\theta}}$- constacyclic code over $\mathtt{R}_{_{\theta}}$ if and only if
			\begin{itemize}
				\item [$(i)$] $h_{_{ii}}(z)$ for $1\leq i \leq 4,$ are all self reciprocal polynomials,
				\item[$(ii)$] $h_{_{44}}(z)\big\vert z^{\alpha}h_{_{34}}^{*}(z)-h_{_{34}}(z),$ where $\alpha=$ deg $h_{_{33}}(z)-$ deg $h_{_{34}}(z) > 0,$
				\item[$(iii)$] $2(z^{\beta}h_{_{12}}^{*}(z)-h_{_{12}}(z)) + k_{_{\theta}}(z^{\gamma} \mathsf{h}_{_{1}}^{*}(z)-\mathsf{h}_{_{1}}(z)) \in \mathtt{C}_{_{\theta}},$ where $\beta =$ deg $h_{_{11}}(z)-$ deg $h_{_{12}}(z) > 0$ and $\gamma =$ deg $h_{_{11}}(z)-$ deg $\mathsf{h}_{_{1}}(z) > 0,$
				\item[$(iv)$] $z^{\delta}\mathsf{h}_{_{2}}^{*}(z)-\mathsf{h}_{_{2}}(z)  \in$ Tor$(\mathtt{C}_{_{\theta}}),$ where $\delta =$ deg $h_{_{22}}(z)-$ deg $\mathsf{h}_{_{2}}(z) \geq 0.$ 
			\end{itemize}
		\end{theorem}
		\begin{proof} The proof follows the same line as the proof of Theorem 3.3 \cite{arxiv rev and rev comp}.
	\end{proof}
	
	The following theorem gives necessary and sufficient conditions for a $\beta_{_{\theta}}$- constacyclic code $\mathtt{C}_{_{\theta}}$ of arbitrary length over $\mathtt{R}_{_{\theta}}, \theta \in \{0,1,\nu\}$ to be a reversible code.
	\begin{theorem}\label{consta beta rev}
		Let $\mathtt{C_{_{\theta}}}$ be a $\beta_{_{\theta}}$- constacyclic code over $\mathtt{R_{_{\theta}}}, \theta \in \{0,1,\nu\}$ of arbitrary length $N=n.2^{s},$ where $s$ and $n$ are positive integers and $n$ is odd. Let $\mathtt{C_{_{\theta}}}= \langle l_{_{\theta_{1}}}(z),l_{_{\theta_{2}}}(z)\rangle,$ where $l_{_{\theta_{i}}}(z) \in \dfrac{\mathtt{R_{_{\theta}}}[z]}{\langle z^{N}-\beta_{_{\theta}} \rangle}$ for $ i=1,2$ are as obtained in Theorem \ref{consta from negacyclic}. Then $\mathtt{C_{_{\theta}}}$ is a reversible $\beta_{_{\theta}}$- constacyclic code if and only if $l_{_{\theta_{1}}}^{*}(z) \in \mathtt{C_{_{\theta}}}.$
	\end{theorem}
	\begin{proof} Let $\mathtt{C_{_{\theta}}}$ be a $\beta_{_{\theta}}$- constacyclic code over $\mathtt{R_{_{\theta}}}, \theta \in \{0,1,\nu\}$ generated by $l_{_{\theta_{1}}}(z)$ and $l_{_{\theta_{2}}}(z).$
		From Lemma \ref{jasbir result on consta}, $\mathtt{C_{_{\theta}}}$ is a reversible $\beta_{_{\theta}}$- constacyclic code if and only if $l_{_{\theta_{1}}}^{*}(z)$ and $l_{_{\theta_{2}}}^{*}(z) \in \mathtt{C_{_{\theta}}}.$ Consider,
		\begin{equation*}
			l_{_{\theta_{2}}}^{*}(z)=\big((z^{n}-1)^{t_{1}}\big)^{*}=(-1)^{t_{1}}(z^{n}-1)^{t_{1}}=(-1)^{t_{1}}l_{_{\theta_{2}}}(z) \in \mathtt{C_{_{\theta}}}
		\end{equation*}  
		It implies that, $\mathtt{C_{_{\theta}}}$ is a reversible $\beta_{_{\theta}}$- constacyclic code if and only if $l_{_{\theta_{1}}}^{*}(z) \in \mathtt{C_{_{\theta}}}.$
	\end{proof}
	
	Some illustrations of Theorem \ref{const from cyc theorem NASC rev} and Theorem \ref{consta beta rev} are given below.
	\begin{example}
		Let $\mathtt{C}_{_{\theta}} = \langle \nu(z+1) \rangle$ be a $(1+\nu)$- constacyclic code of length $3$ over the ring $\mathtt{R}_{_{\theta}}$ for $\theta =0.$ Here $h_{_{11}}(z)=0, h_{_{22}}(z)=0, h_{_{33}}(z)=z+1, h_{_{44}}(z)=0, h_{_{12}}(z)=0, h_{_{34}}(z)=0, \mathsf{h}_{_{1}}(z)=0$ and $\mathsf{h}_{_{2}}(z)=0.$ Clearly, $h_{_{33}}^{*}(z)=z+1=h_{_{33}}(z).$ Hence, $\mathtt{C_{_{\theta}}}$ is a reversible constacyclic code.
	\end{example}
	\begin{example}
		Let $\mathtt{C}_{_{\theta}} = \langle 1+z^{4}+2(1)+\nu(z+1), 2(1+z^{2})+\nu(z-1), \nu(1+z^{4})+2\nu (z-1), 2\nu(z^2+1) \rangle$ be a $(1+2\nu)$- constacyclic code of length $8$ over the ring $\mathtt{R}_{_{\theta}}$ for $\theta =0.$ Here $h_{_{11}}(z)=z^{4}+1, h_{_{22}}(z)=z^{2}+1, h_{_{33}}(z)=1+z^{4}, h_{_{44}}(z)=1+z^{2}, h_{_{12}}(z)=1, h_{_{34}}(z)=z-1, \mathsf{h}_{_{1}}(z)=1+z$ and $\mathsf{h}_{_{2}}(z)=z-1.$ Since $\delta =1$ and $\mathsf{h}_{_{2}}^{*}(z)=z-1,$ it implies that $z\mathsf{h}_{_{2}}^{*}(z)-\mathsf{h}_{_{2}}(z)  \notin$ Tor$(\mathtt{C}_{_{\theta}})$ which violates conditon $(iv)$ of Theorem \ref{const from cyc theorem NASC rev}. Hence, $\mathtt{C_{_{\theta}}}$ is not a reversible constacyclic code.
	\end{example}
	
	\begin{example}
		Let $\mathtt{C_{_{\theta}}}=\langle z-1+\nu,z-1\rangle$ be a $(3+2\nu)$- constacyclic code of length $4$ over the ring $\mathtt{R_{_{\theta}}}$ for $\theta =0.$ Here $n=1,s=2$ and $t_{1}=t_{2}=1.$ Then we have, $l_{_{\theta_{1}}}^{*}(z)=(z-1)^{*}+\nu z=-1(1-z)+\nu z=z(z-1+\nu)+(1-z)(z-1) \in \mathtt{C_{_{\theta}}}.$ Hence, $\mathtt{C_{_{\theta}}}$ is a reversible constacyclic code. 
	\end{example}
	\begin{example}
		Let $\mathtt{C_{_{\theta}}}=\langle (z^{3}-1)^{2}+(1+\nu)z^{2}+2(1+\nu)z,z^{3}-1\rangle$ be a $(1+2\nu)$- constacyclic code of length $3$ over the ring $\mathtt{R_{_{\theta}}}$ for $\theta =1.$ Here $n=3,s=0$ and $t_{1}=t_{2}=1$ Then we have, $l_{_{\theta_{1}}}^{*}(z)=\big((z^{3}-1)^{2}\big)^{*}+(1+\nu)z^{4}(z^{2})^{*}+2(1+\nu)z^{5}z^{*}= (z^{3}-1)^{2}+(1+\nu)z^{4}+2(1+\nu)z^{5}\notin \mathtt{C_{_{\theta}}}.$ Hence, $\mathtt{C_{_{\theta}}}$ is not a reversible constacyclic code. 
	\end{example}
	
	\section{Conclusion and Future Scope}
	\label{}
	In this paper, a unique form of generators of a constacyclic code of arbitrary length over a non-chain ring of the type $\mathtt{R_{_{\theta}}}=Z_{4}+\nu Z_{4}, \nu^{2}=\theta \in Z_{4}+\nu Z_{4}$ has been obtained. Further rank and cardinality of a constacyclic code of arbitrary length over a non-chain ring of the type $\mathtt{R_{_{\theta}}}$ have been determined. In addition, necessary and sufficient conditions for a constacyclic code of arbitrary length over a non-chain ring of the type $\mathtt{R_{_{\theta}}}$ to be reversible have been determined.	
	
	\section{Acknowledgements}
	The first author would like to thank Council of Scientific and Industrial Research (CSIR) India, for providing fellowship.\\ Grant No. 08/423(0003)/2020-EMR-I.
	
\end{document}